 \documentclass[draft]{article}

\usepackage{amsmath,amsfonts,amsthm,amssymb,amscd,cancel,color}
\usepackage{enumitem}
\usepackage{verbatim}
\usepackage[dvipdfmx]{graphicx}
\usepackage{ulem,color}

\renewcommand{\Re}{\mathop{\rm Re}\nolimits}

\theoremstyle{plain}
\newtheorem{theorem}{Theorem}[section]
\newtheorem{lemma}[theorem]{Lemma}
\newtheorem{proposition}[theorem]{Proposition}

\theoremstyle{definition}

\theoremstyle{remark}
\newtheorem{remark}[theorem]{Remark}

\newtheorem{claim}[theorem]{Claim}

\newcommand{\R}{{\mathbb R}}

\def\im{{\rm i}}

\newcommand{\C}{\mathbb{C}}

\newcommand{\sech}{\mathrm{sech}}

\def\({\left(}
\def\){\right)}
\def\<{\left\langle}
\def\>{\right\rangle}
\newcommand{\rad}{{\mathrm{rad}}}

\newcommand{\sechs}{\mathrm{sech}^2}
\newcommand{\phic}{\phi_3}
\newcommand{\sqt}{\sqrt{2}}

\numberwithin{equation}{section}

\begin{document}

\title{A note on the  Fermi Golden Rule  constant for the pure power NLS}

\author{Scipio Cuccagna, Masaya Maeda}
\maketitle

\begin{abstract}
We provide a detailed proof that the Nonlinear Fermi Golden Rule coefficient that appears in our recent  proof of the asymptotic stability of ground states for the pure power Nonlinear Schr\"odinger equations  in $\R$ with  exponent $0<|p-3|\ll 1$ is nonzero.

\end{abstract}

\section{Introduction}
We consider the  pure power Nonlinear Schr\"odinger equation on the line
\begin{align}\label{eq:nls1}
\im \dot{u}-u''  - |u|^{p-1}u=0,\ u:\R^{1+1}\to\C ,
\end{align}
where $\dot{u}=\partial_t u$ and $u'=\partial_x u$.    In this paper  we only consider $p$ near $3$.
It is well known that Equation   \eqref{eq:nls1} has standing waves, solutions with the form $u(t,x)=e^{\im \omega t}\phi_{\omega}(x)$. They are  obtained from  $\phi _p[\omega] (x)=\omega ^{\frac 1{p-1}} \phi  _p (\sqrt{\omega }x) $  with the explicit formula
\begin{align} &
    \phi _p(x) :=  {\(\frac {p+1}2 \)^{\frac 1{p-1}}}{
\sech   ^{\frac 2{p-1}}\(\frac{p-1}2 x\)}  ,\label{eq:sol}
   \end{align}
 see formula (3.1)  Chang et al. \cite{Chang}.
For $p=3$ we have
\begin{align*}
\phic=\sqrt{2}\sech(x).
\end{align*}
It is well known that the \textit{linearization} of \eqref{eq:nls1} at $\phi _p[\omega]$,  is given by
\begin{align}\label{eq:lineariz2} \partial _t   \begin{pmatrix}
 r_1 \\ r_2
 \end{pmatrix}      &=  \mathcal{L}_{p\omega }  \begin{pmatrix}
 r_1 \\ r_2
 \end{pmatrix}  \text{  with }    \mathcal{L}_{p\omega }  := \begin{pmatrix}
0 & L_{-, p,\omega} \\ -L_{+, p,\omega } & 0
\end{pmatrix} ,
    \end{align}
where  \begin{align*}
  L _{+, p,\omega}:=- \partial _x^2   +\omega -p \phi _p  ^{p-1}[\omega]   \text{ and }     L _{-, p,\omega}:=- \partial _x^2   +\omega-  \phi _ p^{p-1}   [\omega].
\end{align*}
The essential spectrum is $\sigma _e (\mathcal{L}_{p\omega } ) =(-\im \infty , -\im \omega ]\cup [\im \omega , \im \infty )$   and $0\in \sigma _p(\mathcal{L}_{p\omega } )$. For $p=3$ there  is no other spectrum.
 Notice that the operator in \eqref{eq:lineariz2} is obtained by a simple scaling from
\begin{align*}
\mathcal{L}_p=\begin{pmatrix}
0 & L_{+,p}\\
-L_{-,p} & 0
\end{pmatrix},
\end{align*}
where
\begin{align*}
L_{-,p}&=-\partial_x^2+1-p\phi_p^{p-1}  \text{ and }
L_{+,p} =-\partial_x^2+1-\phi_p^{p-1}.
\end{align*}
  For  $0<|p-3|\ll 1$,
Coles and Gustafson \cite{coles} proved   $ \mathcal{L}_{p\omega } $ has exactly one eigenvalue of the form  $\im \lambda (p, \omega )  $ near $\im \omega$  where with $ 0<  \lambda (p, \omega ) =\omega \lambda  (p,1) <\omega  $ with  $\dim \ker (\mathcal{L}_{p\omega } -\im \lambda (p,\omega))=1$, thus giving a partial  rigorous confirmation of the numerical results of Chang et al. \cite{Chang}.  For   $\xi _p [\omega]\in H^1(\R, \C^2) $ a generator of $\ker (\mathcal{L}_{p\omega } -\im \lambda (p,\omega))$   in  \cite{CM24D1} we stated the following    result involving only radial functions in $H^1_\rad (\R , \C) $.

\begin{theorem} \label{thm:asstab}
There exists $p_1<3<p_2$ s.t.\ for any $p\in (p_1,p_2)\setminus\{3\}$ and any $\omega _0 >0$, any $a>0$ and any $\epsilon >0$   there exists a  $\delta >0$ such that for any  initial value    $u_0\in  D _{H^1_\rad (\R ) } ({\phi}_{\omega  _0},\delta  ) $
there exist functions $ \vartheta , \omega \in C^1 \( \R , \R \)$,  $ z \in C^1 \( \R , \C  \)$ and $\omega_+>0$    s.t.\  the  solution  of \eqref{eq:nls1} with initial datum  $u_0$  can be written as
 \begin{align} \label{eq:asstab1}
    & u(t)= e^{\im \vartheta (t)} \(   \phi _{\omega (t)} + z(t) \xi _p [ \omega(t) ]+ \overline{z}(t) \overline{\xi} _p [ \omega(t) ] +\eta (t)\) \text{   with}
 \\&  \label{eq:asstab2}   \int _{\R } \|  e^{- a\< x\>}   \eta (t ) \| _{H^1(\R )}^2 dt <  \epsilon  \text{  where }\< x\>:=\sqrt{1+x^2} ,\\&  \label{eq:asstab20}
   \lim _{t\to \infty  }  \|  e^{- a\< x\>}   \eta (t ) \| _{L^2(\R )}     =0 ,  \\&  \label{eq:asstab2}
   \lim _{t\to \infty}z(t)=0      ,\\&
%\end{align}
%Furthermore, for $p>3$  there are constants $\omega _\pm >$ s.t. we have
%\begin{align}
\lim _{t\to \infty}\omega (t)= \omega _+ .\label{eq:asstab3}
\end{align}

\end{theorem}

\begin{remark}\label{rem:asstab3--}  Notice that the fact that the $H^1(\R )$ norm of $\eta$ is uniformly bounded for all times, guaranteed by classically known  orbital stability  results  of Cazenave and Lions \cite{cazli},  Shatah  \cite{shatah} and Weinstein  \cite{W1}, and
Theorem \ref{thm:asstab} imply  the following local in space asymptotic  convergence up to the phase,
\begin{align}\label{eq:rem:asstab3--1}
  \lim _{t\to +\infty }  u(t) e^{-\im \vartheta (t)}=  \phi  _{\omega _+}  \text{  in }L^\infty _{\text{loc}}(\R ) .
\end{align}
\end{remark}
A result similar to Theorem \ref{thm:asstab} was obtained completely independently,   unbeknownst to the two sets of authors during the process of writing the respective papers  and posted on Arxiv the same day,
  with a different proof by Rialland \cite{Rialland2}, which remains closer than us to the  theory in Martel \cite{Martel2}, which in turn has a  similar result but for a cubic quintic NLS.

We set
\begin{align*}
   \xi _p [ \omega  ] =\(  \xi _{p,1} [ \omega  ] ,  \xi _{p,2} [ \omega  ]     \) ^\intercal
\end{align*}
and let furthermore  $g ^{(\omega)}_p =\( g_{p,1} ^{(\omega)}  , g_{p,2}  ^{(\omega)} \) ^\intercal \in L^\infty \(   \R , \C^2\) $  be an nonzero solution of
\begin{align}
  \label{eq:eqsatg2} \mathcal{L}_{p\omega}g ^{(\omega)}_p=2\im \lambda (p, \omega) g ^{(\omega)}_p,
\end{align}
see     Buslaev and Perelman \cite{BP1} and  Krieger and Schlag \cite{KrSch}.    Notice that  if $g _p = g ^{(1)} _p$,     then $g ^{(\omega)} (x) = g_p\( \sqrt{\omega}x \)$      solves
\eqref{eq:eqsatg2}.  Similarly  if $\xi _p = \xi   _p [1] $,     then $\xi   _p [\omega]  (x) = \xi_p\( \sqrt{\omega}x \)$      is a generator of  $\ker (\mathcal{L}_{p\omega } -\im \lambda (p,\omega))$. Here, for later reference, we record that
\begin{align}
  \label{eq:eqsatg1} \mathcal{L}_{p }\xi _p = \im \lambda (p, 1) \xi _p.
\end{align}
  In the special case $p=3 $ all the solutions to \eqref{eq:eqsatg1} have been known explicitly    since
Kaup \cite{kaup}.
  Using now  the notation
\begin{align*}
  \<f,g\> := \Re \int _\R f(x) \overline{g(x)} dx ,
\end{align*}
we consider  the constant
\begin{equation}\label{eq:fgrgamma}
   \gamma ( p, \omega )  :=     \<    \phi _p ^{p-2} [ \omega  ]   \( p  \xi_{p,1} ^2 [ \omega  ] +  \xi_{p,2}  ^2  [ \omega  ] \) , g _{p,1} ^{(\omega)} \> +2 \<    \phi _p ^{p-2}  [ \omega  ]   \xi_{p,1} [ \omega  ]  \xi_{p,2} [ \omega  ]  , g_{p,2}  ^{(\omega)}\>  ,
\end{equation}
which plays an important  role in the proof of Theorem \ref{thm:asstab}.
Let $\gamma ( p  ):=\gamma ( p, 1 )$.
In \cite{CM24D1} we proved that it is possible to choose  $g_p$ and $\xi _p$  so that    there is a constant $\Gamma \in \R$  such that
\begin{align}\label{eq:espgamma}
   \gamma ( p  ) =   (p-3) \Gamma + O\(  (p-3)^2   \).
\end{align}
The purpose of the present paper is to prove   the following.
\begin{theorem}
  \label{thm:fgr} We have     \begin{align} \label{eq:fgr1}
\Gamma=\frac{\pi}{\sqrt{2}\cosh(\pi/2)}  .
\end{align}
\end{theorem}

This implies automatically that the corresponding constant  \eqref{eq:fgrgamma} is nonzero for any $\omega >0$ and for $0<|p-3|\ll 1$.

\noindent The proof of Theorem \ref{thm:asstab} in \cite{CM24D1} utilizes in an essential way that $ \gamma ( p  ) \neq 0$ for   $0<|p-3|\ll 1$ and gives  only a sketch of Theorem \ref{thm:fgr}. A similar sketch is given in Rialland \cite{Rialland2}.
Both  \cite{CM24D1,Rialland2} follow closely  Martel \cite{Martel2} for their proof of their own versions of  Theorem \ref{thm:fgr} and all three papers omit to write a full proof, on account of  this would  add many  pages of very elementary  computations  to papers devoted to solving longstanding problems  in the context of classical models like Equation \eqref{eq:nls1}.  Our aim in this paper is to fill in the details of the proof of Theorem \ref{thm:fgr},  following the succinct description of the analogous computation in   Martel \cite{Martel2}.  While the computations added here and missing in  \cite{CM24D1} are elementary, we think  that writing the complete proof    might result in  a useful reference   for  fellow researches.   Besides, the computation of the constant
$ \Gamma$ raises at least one interesting question.
Indeed, it is rather surprising that both  for the pure power equation \eqref{eq:nls1} and for the cubic--quintic equation in Martel \cite{Martel2}, the constant   $ \Gamma$ has a very simple form. We ignore  why this is the case, but certainly this merits further investigation. Rialland \cite{Rialland2} has also some interesting numerical computations that show that $ \gamma ( p  ) \neq 0$ further away from $p=3$.

   We explain now, briefly, why it is important that these constants be different from 0.  In the study of stability problems like  in Theorem \ref {thm:asstab} it is quite natural to linearize
Hamiltonian systems \eqref{eq:nls1} at ground states. The linearized operators have both discrete and continuous modes and it is quite natural to use these  modes to study also the nonlinear systems. At a linear level the different modes are uncoupled and so the discrete modes tend to describe harmonic oscillators. However there is mechanism, first introduced by Sigal  \cite{sigal}, the so called nonlinear Fermi Golden Rure (from now on FGR), which explains why,  by nonlinear interaction, the discrete modes, specifically related to excited states of the linearization,  lose energy. The discrete modes are like a finite dimensional Hamiltonian system embedded in a sea of  radiation represented by the continuous modes, with energy   spilling out of the finite dimensional system and then, essentially by linear dispersion,    scattered to infinity. Sigal's idea was developed in relatively simple contexts by Buslaev and Perelman \cite{BP2}, see also Buslaev and Sulem \cite{BS}, and by Soffer  and Weinstein \cite{SW3}. Later \cite{Cu1} described the FGR for general spectral configurations, illustrating the deep relation between the FGR and the Hamiltonian structure of the systems. The FGR involves the fact that certain coefficients of the mixed system of discrete and continuous modes, have nonzero restriction on appropriate sphere of the \textit{distorted} phase space associated to the linearization \eqref{eq:lineariz2}. In dimension 1  the spheres reduce  to pair of points and the FGR reduces to having certain  single frequencies different from 0, and that is exactly the meaning of  the condition $ \gamma ( p, \omega ) \neq 0$.

\section{A first formula for the constant  $ \Gamma$} \label{sec:rew}

The main purpose of this section is to prove Proposition \ref{lem:1st}.  Before stating and proving it, we   review briefly some results on $ \gamma ( p  )$ proved in  \cite{CM24D1}. First of all we proved the following.

\begin{lemma}
	\label{lem:existeig} There exists a small $\delta _1>0$ and a function   $\alpha \in C^\infty (D _{\R} (3,\delta _1) , \R )$ such that
\begin{align}\nonumber
\lambda(p,1) &=1-\alpha  ^2(p) \text{ with }
 \\\label{eq:expansalpha}
		\alpha (p) &=   (p-3) ^2 c_0   + O\( (p-3) ^3    \) \text{ where } c_0 =   2^{-2} + 2^{-5}2^{-\frac{1}{2}}\< \phi_3^2,  { T}        \>,
	\end{align}
where
\begin{equation}\label{eq:defboldT}
	 {T}:=  e^{- \sqrt{2}  |\cdot | }   * \sech ^2 (x) .
\end{equation}
	
\end{lemma}

\qed

Notice that the above was a different proof, simpler thanks to the  Darboux transformation due to Martel \cite{Martel2}, of the result on the eigenvalue by Coles and  Gustafson  \cite{coles}. The specific value of the constant $c_0>0$ does not play a role here and \cite{CM24D1}.

\noindent Next, in     \cite{CM24D1} we proved the following.

\begin{lemma}	\label{lem:expv} There exists an open interval $\mathcal{I}$ containing 3 and  for each $p\in \mathcal{I}$
	it is possible to choose $ \xi_p=(\xi_{p,1},\xi_{p,2})^\intercal$ so that \eqref{eq:eqsatg1} holds,
we have
\begin{align*}
 \xi_3=\(1-\phi_3^2,\im\) ^\intercal
\end{align*}
and
	\begin{align}
		\xi _{p,1}&=1-\phi_3^2+(p-3)R_1+(p-3)^2\tilde{\xi} _{p,1},\label{eq:expandxip1}\\
		\xi _{p,2} &=\im\(1+(p-3)R_2+(p-3)^2\tilde{\xi} _{p,2}\),\label{eq:expandxip2}
	\end{align}
	where
	\begin{align*}
		R_1&=-x\phi_3\phi_3'-\frac{1}{4\sqrt{2}}(3-\phi_3^2) {T}-\frac{\phi_3'}{2\sqrt{2}\phi_3} { T}' \text{ and }\\
		R_2&=\frac{1}{2}\phi_3^2+\frac{3}{4\sqrt{2}} {T}+\frac{\phi_3'}{2\sqrt{2}\phi_3} { T}'
	\end{align*}
	and where furthermore, for any $k\ge 0$ there exists a constant $C_k$ such that
	\begin{align} \label{eq:exptildev}
		& |  \widetilde{\xi} _{p,j} ^{(k)} (x) |\le  C_k  \< x\>  ^{3} \text{ for all } x\in \R \text{ and all $p\in \mathcal{I}$}. \end{align} Here $R_j$  and $\tilde{\xi} _{p,j}$   for $j=1,2$ are real valued.

\end{lemma}
\qed

Next, we  set
\begin{align}
  \label{eq:defg3}g_3=  \( h _{3,1},\im  h _{3,2}  \) ^\intercal \text{ with } \( h _{3,1},   h _{3,2}  \) ^\intercal
  =
  \(\frac{1}{2}\phi_3^2\cos(x)+\frac{\phi_3'}{\phi_3}\sin(x), \frac{\phi_3'}{\phi_3}\sin(x)\) ^\intercal .
\end{align}
Then  in \cite{CM24D1} we stated without   proof the following    variant of Lemma 19 of Martel \cite{Martel2}, which can be proved similarly.
\begin{lemma}There exists an open interval $\mathcal{I}$ containing 3 and  for each $p\in \mathcal{I}$
	it is possible to choose  $g_p$ so that
\begin{align*}
\|g_p-\(\frac{1}{2}\phi_3^2\cos(\tau x)+\frac{\phi_3'}{\phi_3}\sin(\tau x),\im \frac{\phi_3'}{\phi_3}\sin(\tau x)\)\|_{L^\infty}\lesssim |p-3|,
\end{align*}
where $\tau=\sqrt{1-\lambda(p,1)^2}$  and  \begin{align*}
   |\partial _x ^ k g_p(x) | \le C _k  \text{  and }  |\partial _x ^ k\partial _p  g_p(x) |\le C _k( 1+|x|) \text{  for all  } x\in \R \text{ and all $p\in \mathcal{I}$}.
\end{align*}
\end{lemma}
\qed

An elementary   differentiation yields
\begin{align*}
E&:=\left.\partial_p\right|_{p=3}\phi_p=\frac{1}{2}\phi_3\(\frac{1}{4}-\log \phi_3\) +\frac{1}{2}x \phi' _3,\\
F&:=\left.\partial_p\right|_{p=3}\phi_p^{p-2}=E+\phi_3\log\phi_3.
\end{align*}
Recalling \begin{align}\label{eq:gammap}
\gamma(p )=\<\phi_p^{p-2}(p\xi_{p,1}^2+\xi_{p,2}^2),g_{p,1}\>+2\<\phi_p^{p-2}\xi_{p,1}\xi_{p,2},g_{p,2}\> ,
\end{align}
we   set, following the notation and the argument  in Martel \cite{Martel2},
\begin{align*}
G_{p,1}:=\phi_p^{p-2}(p\xi_{p,1}^2+\xi_{p,2}^2),\\
G_{p,2}:=-2\im \phi_p^{p-2}\xi_{p,1}\xi_{p,2}.
\end{align*}
Then by elementary computations
\begin{align*}
G_{p,1}&=\phi_3\(3(1-\phi_3^2)^2-1\)+(p-3)\Delta_1+\tilde{\Delta}_1,\\
\frac{1}{2}G_{p,2}&=\phi_3(1-\phi_3^2)+(p-3)\Delta_2+\tilde{\Delta}_2,
\end{align*}
 where
\begin{align*}
\Delta_1&=F\(3(1-\phi_3^2)^2-1\) + \phi_3(1-\phi_3^2)^2+6\phi_3(1-\phi_3^2)R_1-2\phi_3R_2,\\
\Delta_2&=F(1-\phi_3^2)+\phi_3R_1+\phi_3(1-\phi_3^2)R_2
\end{align*}
and $\tilde{\Delta}_1, \tilde{\Delta}_2$ are   remainder terms of $(p-3)^2$ order.
The following formula is proved in \cite{CM24D1},
\begin{align}\nonumber
 \Gamma&= \<\Delta_1,h_{3,1}\> +2\<\Delta_2,h_{3,2}\> + \frac{1}{2}\< \(6x\tanh \sechs -\frac{7}{2}\sechs\)G_{3,2},h_{3,1} \> -2\<E,h_{3,1}\> \\& =:\gamma_1+\gamma_2+\gamma_3+\gamma_4,  \label{eq:decgamma1}
\end{align}
where  $\sech  =\sech (x) $, $\log =\log (x)$, $\tanh = \tanh (x)$ etc. With this notation,
denoting similarly  $\log\circ \sech=\log \(  \sech (x) \)$, $T=T(x)$,
we   consider the following constants, similar to analogous ones introduced in  Martel \cite{Martel2},
 \begin{align*}
p_k&=\int\sech^k\cos,\\
q_k&=\int \sech^k  \log\circ \sech   \cos \\
r_k&=\int \sech^k   {T}\cos\\
s_k&=\int \sech^k {T} \tanh \sin\\
a_k&=\int x \sech^k \tanh \cos\\
b_k&=\int \sech^k\tanh \sin\\
c_k&=\int \sech^k \log\circ \sech  \tanh\sin\\
d_k&=\int x \sech^k  \sin\\
e_k&=\int \sech^k \tanh   {T}'\cos\\
f_k&=\int \sech^k  {T}'\sin    .
\end{align*}
The first step in our computations is the following.
\begin{proposition}\label{lem:1st} The following formula holds,
 \small   \begin{align*}
\Gamma=&\sqt\(\(\frac{1}{2}\log 2 + \frac{17}{4}\)p_1 +  \(-13\log 2 -71\)p_3+ \(28\log 2 +\frac{311}{2}\)p_5 + \(-15\log 2 -\frac{173}{2}\)p_7\)\\&
+\sqt\(3q_1 -28q_3 +56 q_5 -30q_7 - a_1 + 93 a_3-336 a_5 +252 a_7\)\\&
+4r_1 -2 r_3 -30 r_5 +30r_7 +2s_1 +  18 s_3 -20 s_5.
\end{align*} \normalsize
\end{proposition}

In order to prove the proposition we will examine separately
each term in \eqref{eq:decgamma1}.
\begin{lemma}\label{lem:gamma1} We have \small
  \begin{align*}
\gamma_1 =&\sqt\(  \(-3\log2-\frac{15}{2}\)p_5+\(3\log2+\frac{11}{2}\) p_7\)\\& +
\sqt\(  q_3 -6q_5+6q_7 -a_3 +  18 a_5 - 30 a_7\)\\& +\sqt\(  \( 3\log 2+\frac{15}{2}\)b_3 -\( 3\log 2 +\frac{11}{2}\)b_5\)\\&
+\sqt\( -c_1 +6 c_3 -6c_5 +d_1 -19 d_3+ 48 d_5-30 d_7\)\\&-6r_3+12r_5-6r_7 +6s_1-12s_3+6s_5  +4e_3-6e_5 -4 f_1  + 10f_3-6f_5.
\end{align*}\normalsize
\end{lemma}
\proof
We write \begin{align*}
\gamma_1&=\<\Delta_1,h_{3,1}\>\\&=\<F(3\xi_{3,1}^2-\xi_{3,2}^2)+\phi_3 \xi_{3,1}^2+6\phi_3\xi_{3,1}R_1-2\phi_3 \xi_{3,2}R_2,h_{3,1}\>\\&
=3\<F\xi_{3,1}^2,h_{3,1}\>-\<F\xi_{3,2}^2,h_{3,1}\>+\<\phi_3\xi_{3,1}^2,h_{3,1}\>+6\<\phi_3\xi_{3,1}R_1,h_{3,1}\>-2\<\phi_3\xi_{3,2}R_2,
h_{3,1}\>\\&
=:\gamma_{11}+\gamma_{12}+\gamma_{13}+\gamma_{14}+\gamma_{15}.
\end{align*}
\begin{claim}
   We have
   \begin{align*}
      \gamma_{11}&=  \gamma_{111}+\gamma_{112}
   \end{align*}
   with
   \begin{align}\nonumber  \gamma_{111}=&
      -\frac{3}{\sqt}(2\log 2+1)p_5+\frac{3}{\sqt}(2\log 2+1)p_7
+\frac{3}{\sqt}q_3-\frac{12}{\sqt}q_5+\frac{12}{\sqt}q_7
\\& -\frac{3}{\sqt}a_3+\frac{12}{\sqt}a_5-\frac{12}{\sqt}a_7 \label{eq:gamma111}
   \end{align}
   and
   \begin{align} \nonumber \gamma_{112}&= \frac{3}{\sqt}(2\log 2+1)b_3-\frac{3}{\sqt}(2\log 2+1)b_5
-\frac{3}{\sqt}c_1+\frac{12}{\sqt}c_3-\frac{12}{\sqt}c_5\\&\quad
+\frac{3}{\sqt}(d_1-d_3)-\frac{12}{\sqt}(d_3-d_5)+\frac{12}{\sqt}(d_5-d_7) . \label{eq:gamma112}
   \end{align}
\end{claim}
\proof  We have
  \small \begin{align*}
\gamma_{11} &=3\<F\xi_{3,1}^2,h_{3,1}\>\\
&=3\<\(\frac{1}{4\sqt}(2\log 2+1)\sech +\frac{1}{\sqt}\sech \log\circ\sech -\frac{1}{\sqt}\sech\cdot x\tanh\)(1-2\sechs)^2,\sechs\cos-\tanh\sin\>\\&
=3\<\(\frac{1}{4\sqt}(2\log 2+1)\sech +\frac{1}{\sqt}\sech \log\circ\sech -\frac{1}{\sqt}\sech\cdot x\tanh\)(1-2\sechs)^2,\sechs\cos\>\\&\quad
-3\<\(\frac{1}{4\sqt}(2\log 2+1)\sech +\frac{1}{\sqt}\sech \log\circ \sech -\frac{1}{\sqt}\sech\cdot x\tanh\)(1-2\sechs)^2,\tanh\sin\>\\&
=:\gamma_{111}+\gamma_{112}.
\end{align*}\normalsize
where we warn the reader that $\sech\cdot x$, here and in analogous expressions below,  stands for the product of the functions $\sinh (x)$ and  $  x$.
We have
\small\begin{align*}
\gamma_{111} &=\frac{3}{\sqrt{2}}\<\(\frac{1}{4 }(2\log 2+1)\sech + \sech \log\circ \sech - \sech\cdot x\tanh\)(1-2\sechs)^2,\sechs\cos\>\\&
=3\<\(\boxed{\frac{1}{4\sqt}(2\log 2+1)\sech} +\frac{1}{\sqt}\sech \log\circ  \sech -\frac{1}{\sqt}\sech\cdot x\tanh\),\sechs\cos\>\\&
\quad -12\<\(\frac{1}{4\sqt}(2\log 2+1)\sech +\frac{1}{\sqt}\sech \log\circ  \sech -\frac{1}{\sqt}\sech\cdot x\tanh\)\sechs,\sechs\cos\>
\\&\quad
+12\<\(\frac{1}{4\sqt}(2\log 2+1)\sech +\frac{1}{\sqt}\sech \log\circ \sech -\frac{1}{\sqt}\sech\cdot x\tanh\)\sech^4,\sechs\cos\>\\&
=-\frac{3}{\sqt}(2\log 2+1)\int \sech^5\cos+\frac{3}{\sqt}(2\log 2+1)\int \sech^7\cos\\&\quad
+\frac{3}{\sqt}\int \sech^3   \log\circ \sech \cos-\frac{12}{\sqt}\int \sech^5   \log\circ \sech \cos+\frac{12}{\sqt}\int \sech^7 \cdot \log\circ\sech \cos\\&\quad
-\frac{3}{\sqt}\int \sech^3\cdot x\tanh \cos+\frac{12}{\sqt}\int \sech^5\cdot x\tanh \cos-\frac{12}{\sqt}\int \sech^7\cdot x\tanh \cos\\&
\end{align*}
\normalsize
where      we get \eqref{eq:gamma111}  if we redefine $\gamma_{111}$  as the last term in the formula,
where, by an abuse of notation, in the last term we did not rewrite the   boxed part, for reasons of space and  because    by $\< \phi _3, h_{3,1}\> =0$ it cancels when added to the boxed part in the next formula. The next formula is, applying the same convention,
 \small
\begin{align*}
\gamma_{112}&=-3\<\(\frac{1}{4\sqt}(2\log 2+1)\sech +\frac{1}{\sqt}\sech \log\circ \sech -\frac{1}{\sqt}\sech\cdot x\tanh\)(1-2\sechs)^2,\tanh\sin\>\\&
=-3\<\(\boxed{\frac{1}{4\sqt}(2\log 2+1)\sech}  +\frac{1}{\sqt}\sech \log\circ  \sech -\frac{1}{\sqt}\sech\cdot x\tanh\),\tanh\sin\>\\&
\quad +12\<\(\frac{1}{4\sqt}(2\log 2+1)\sech +\frac{1}{\sqt}\sech \log\circ  \sech -\frac{1}{\sqt}\sech\cdot x\tanh\)\sechs,\tanh\sin\>
\\&\quad
-12\<\(\frac{1}{4\sqt}(2\log 2+1)\sech +\frac{1}{\sqt}\sech \log\circ \sech -\frac{1}{\sqt}\sech\cdot x\tanh\)\sech^4,\tanh\sin\>\\&
=\frac{3}{\sqt}(2\log 2+1)\int \sech^3\tanh\sin-\frac{3}{\sqt}(2\log 2+1)\int \sech^5\tanh\sin\\&\quad
-\frac{3}{\sqt}\int \sech   \log\circ  \sech \tanh\sin+\frac{12}{\sqt}\int \sech^3   \log\circ  \sech \tanh\sin-\frac{12}{\sqt}\int \sech^5   \log\circ  \sech \tanh\sin\\&\quad
+\frac{3}{\sqt}\int \sech\cdot x(1-\sechs)\sin-\frac{12}{\sqt}\int \sech^3\cdot x(1-\sechs)\sin+\frac{12}{\sqt}\int \sech^5\cdot x(1-\sechs)\sin,
\end{align*}
\normalsize where  we get \eqref{eq:gamma112}.
\qed

\begin{claim}
   We have
   \begin{align*}
     \gamma_{12}&
=\gamma_{121}+\gamma_{122}
   \end{align*}
 with
   \begin{align}   \gamma_{121}&=-\frac{1}{\sqt}q_3
+\frac{1}{\sqt}a_3,
      \label{eq:gamma121}
   \\   \gamma_{122}&=\frac{1}{\sqrt{2}}c_1
-\frac{1}{\sqt}(d_1-d_3)   . \label{eq:gamma122}
   \end{align}

\end{claim}
\proof We have
\begin{align*}
\gamma_{12}&=-\<F\xi_{3,2}^2,h_{3,1}\>\\&
=-\<\cancel{\frac{1}{4\sqt}(2\log 2+1)\sech}  +\frac{1}{\sqt}\sech \log \circ \sech -\frac{1}{\sqt}\sech\cdot x\tanh,\sechs\cos-\tanh\sin\>\\&
=-\< \frac{1}{\sqt}\sech \log \circ  \sech -\frac{1}{\sqt}\sech\cdot x\tanh,\sechs\cos\>\\&\quad
+\<\frac{1}{\sqt}\sech \log \circ  \sech -\frac{1}{\sqt}\sech\cdot x\tanh,\tanh\sin\>\\&
=:\gamma_{121}+\gamma_{122},
\end{align*}
where the canceled term is null by    $\< \phi _3, h_{3,1}\> =0$. Then from
\begin{align*}
\gamma_{121}&=-\<\frac{1}{\sqt}\sech \log\circ \sech -\frac{1}{\sqt}\sech\cdot x\tanh,\sechs\cos\>\\&
=
-\frac{1}{\sqt}\int \sech^3 \log\circ \sech \cos
+\frac{1}{\sqt}\int\sech^3\cdot x\tanh\cos
\end{align*}
and
\begin{align*}
\gamma_{122}&=\< \frac{1}{\sqt}\sech \log\circ \sech -\frac{1}{\sqt}\sech\cdot x\tanh,\tanh\sin\>\\&
=\frac{1}{\sqt}\int \sech \log\circ \sech  \tanh\sin
-\frac{1}{\sqt}\int \sech(1-\sechs)  x\sin,
\end{align*}
we obtain the desired claim.

\qed

\begin{claim}
   We have
   \begin{align*}
     \gamma_{13}&
=\gamma_{131}+\gamma_{132}
   \end{align*}
 with
   \begin{align}   \gamma_{131}&=-4\sqt p_5+4\sqt p_7,
      \label{eq:gamma131}
   \\   \gamma_{132}&=
4\sqt b_3 -4\sqt b_5   . \label{eq:gamma132}
   \end{align}

\end{claim}
\proof We have
\begin{align*}
\gamma_{13}&=\<\phi_3\xi_{3,1}^2,h_{3,1}\>\\&
=\sqt\<\sech (1-2\sechs)^2,\sechs\cos-\tanh\sin\>\\&
=\sqt\<\cancel{\sech} -4\sech^3+4\sech^5,\sech^2\cos-\tanh\sin\>\\&
=\sqt\<-4\sech^3+4\sech^5,\sech^2\cos\>
-\sqt\<-4\sech^3+4\sech^5,\tanh\sin\>\\&
=:\gamma_{131}+\gamma_{132},
\end{align*}
where the canceled term is null by    $\< \phi _3, h_{3,1}\> =0$.
Then from
\begin{align*}
\gamma_{131}&=\sqt\<-4\sech^3+4\sech^5,\sech^2\cos\>\\&
=-4\sqt\int\sech^5\cos+4\sqt\int\sech^7\cos
\end{align*}
and
\begin{align*}
\gamma_{132}&=-\sqt\<-4\sech^3+4\sech^5,\tanh\sin\>\\&
=4\sqt\int \sech^3\tanh\sin -4\sqt\int \sech^5\tanh\sin,
\end{align*}
we obtain the claim.
 \qed

\begin{claim}
   We have
   \begin{align*}
     \gamma_{14}&
=\gamma_{141}+\gamma_{142}
   \end{align*}
 with
   \begin{align}   \gamma_{141}&= 12\sqt (a_5-2a_7)-\frac{3}{2}(3r_3-8r_5+4r_7)+3(e_3-2e_5)
      \label{eq:gamma141}
   \\   \gamma_{142}&=-12\sqt(d_3-3d_5+2d_7)+\frac{3}{2}\(3s_1-8s_3+4s_5\)-3\(f_1-3f_3+2f_5\)  . \label{eq:gamma142}
   \end{align}

\end{claim}
\proof We have
\small
\begin{align*}
\gamma_{14}&=6\<\phi_3\xi_{3,1}R_1,h_{3,1}\>\\&
=6\sqt\<\sech (1-2\sechs)\(2x\tanh \sech^2-\frac{1}{4\sqt}(3-2\sechs)T+\frac{1}{2\sqt}\tanh T'\),\sechs\cos-\tanh\sin\>\\&
=6\sqt\<\sech (1-2\sechs)\(2x\tanh \sech^2-\frac{1}{4\sqt}(3-2\sechs)T+\frac{1}{2\sqt}\tanh T'\),\sechs\cos\>\\&\quad
-6\sqt\<\sech (1-2\sechs)\(2x\tanh \sech^2-\frac{1}{4\sqt}(3-2\sechs)T+\frac{1}{2\sqt}\tanh T'\),\tanh\sin\>\\&
=:\gamma_{141}+\gamma_{142}.
\end{align*} \normalsize
From
\begin{align*}
\gamma_{141}&=6\sqt\<\sech (1-2\sechs)\(2x\tanh \sech^2-\frac{1}{4\sqt}(3-2\sechs)T+\frac{1}{2\sqt}\tanh T'\),\sechs\cos\>\\&
=12\sqt \int (\sech^5-2\sech^7)x\tanh\cos-\frac{3}{2}\int(3\sech^3-8\sech^5+4\sech^7)T\cos\\&\quad
+3\int(\sech^3-2\sech^5)\tanh T'\cos
\end{align*}
and
\begin{align*}
\gamma_{142}&=-6\sqt\<\sech (1-2\sechs)\(2x\tanh \sech^2-\frac{1}{4\sqt}(3-2\sechs)T+\frac{1}{2\sqt}\tanh T'\),\tanh\sin\>\\&
=-12\sqt\int (\sech^3-3\sech^5+2\sech^7)x\sin+\frac{3}{2}\int (3\sech-8\sech^3+4\sech^5)T\tanh\sin\\&\quad
-3\int\(\sech-3\sech^3+2\sech^5\)T'\sin,
\end{align*}
we obtain the claim.
  \qed

\begin{claim}
   We have
   \begin{align*}
     \gamma_{15}&
=\gamma_{151}+\gamma_{152}
   \end{align*}
 with
   \begin{align}   \gamma_{151}&=  -2\sqt p_5-\frac{3}{2}r_3+e_3,
      \label{eq:gamma151}
   \\   \gamma_{152}&=   2\sqt b_3 +\frac{3}{2} s_1 - f_1+f_3. \label{eq:gamma152}
   \end{align}

\end{claim}
\proof
We have  %\small
\begin{align*}
\gamma_{15} &=-2\<\phi_3\xi _{3,2}R_2,h_{3,1}\>\\&
=-2\sqt\<\sech \(\sechs  +\frac{3}{4\sqt}T-\frac{\tanh}{2\sqt}T'\),\sechs\cos-\tanh\sin\>\\&
=-2\sqt\<\sech \(\sechs  +\frac{3}{4\sqt}T-\frac{\tanh}{2\sqt}T'\),\sechs\cos\>
\\&\quad+2\sqt\<\sech \(\sechs  +\frac{3}{4\sqt}T-\frac{\tanh}{2\sqt}T'\),\tanh\sin\>\\&
=:\gamma_{151}+\gamma_{152}.
\end{align*}
%\normalsize
From
\begin{align*}
\gamma_{151}&=-2\sqt\<\sech \(\sechs  +\frac{3}{4\sqt}T-\frac{\tanh}{2\sqt}T'\),\sechs\cos\>\\&
=-2\sqt\int \sech^5\cos-\frac{3}{2}\int \sech^3 T\cos+\int \sech^3 \tanh T' \cos
\end{align*}
and
\begin{align*}
\gamma_{152}&=2\sqt\<\sech \(\sechs  +\frac{3}{4\sqt}T-\frac{\tanh}{2\sqt}T'\),\tanh\sin\>\\&
=2\sqt \int \sech^3 \tanh \sin +\frac{3}{2}\int \sech T \tanh \sin -\int \sech (1-\sechs)T'\sin ,
\end{align*}
we obtain the claim.
  \qed

From
\begin{align*}
\gamma_1&=\gamma_{111}+\gamma_{112}+\gamma_{121}+
\gamma_{122}+\gamma_{131}+\gamma_{132}+\gamma_{141}+\gamma_{142}+\gamma_{151}+\gamma_{152},
\end{align*}
we obtain \small
\begin{align*}
\gamma_1&
=-\frac{3}{\sqt}(2\log 2+1)p_5+\frac{3}{\sqt}(2\log 2+1)p_7
+\frac{3}{\sqt}q_3-\frac{12}{\sqt}q_5+\frac{12}{\sqt}q_7
-\frac{3}{\sqt}a_3+\frac{12}{\sqt}a_5-\frac{12}{\sqt}a_7\\&\quad
+\frac{3}{\sqt}(2\log 2+1)b_3-\frac{3}{\sqt}(2\log 2+1)b_5
-\frac{3}{\sqt}c_1+\frac{12}{\sqt}c_3-\frac{12}{\sqt}c_5\\&\quad
+\frac{3}{\sqt}(d_1-d_3)-\frac{12}{\sqt}(d_3-d_5)+\frac{12}{\sqt}(d_5-d_7)\\&\quad
\underbrace{
-\frac{1}{\sqt}q_3
+\frac{1}{\sqt}a_3}_{=\gamma_{121}}
+\underbrace{\frac{1}{\sqt}c_1
-\frac{1}{\sqt}(d_1-d_3)}_{=\gamma_{122}}\\&\quad
\underbrace{-4\sqt p_5+4\sqt p_7}_{=\gamma_{131}}+\underbrace{4\sqt b_3 - 4\sqt b_5}_{=\gamma_{132}}+\underbrace{12\sqt (a_5-2a_7)-\frac{3}{2}(3r_3-8r_5+4r_7)+3(e_3-2e_5)}_{=\gamma_{141}}\\&\quad
\underbrace{-12\sqt(d_3-3d_5+2d_7)+\frac{3}{2}\(3s_1-8s_3+4s_5\)-3\(f_1-3f_3+2f_5\)}_{=\gamma_{142}}\\&\quad
\underbrace{-2\sqt p_5-\frac{3}{2}r_3+e_3}_{=\gamma_{151}}+\underbrace{2\sqt b_3 +\frac{3}{2} s_1 - f_1+f_3}_{=\gamma_{152}}
\end{align*}
\normalsize
Summing up corresponding terms and tracking them to make the computations simpler to read, we obtain
\small
\begin{align*}&
\gamma_1 =\sqt\(  \underbrace{\(-\frac{3}{2}(2\log 2+1)-4-2\)}_{=-3\log2-\frac{15}{2}}p_5+\underbrace{\(\frac{3}{2}(2\log 2+1)+4\)}_{=3\log2+\frac{11}{2}}p_7\)\\& +
\sqt\( \underbrace{\(\frac{3}{2}-\frac{1}{2}\)}_{=1}q_3 -6q_5+6q_7 +\underbrace{\(-\frac{3}{2}+\frac{1}{2}\)}_{=-1}a_3 + \underbrace{\(6+12\)}_{=18}a_5+\underbrace{\(-6-24\)}_{=-30}a_7\)\\& +\sqt\( \underbrace{\(\frac{3}{2}(2\log 2+1)+4+2\)}_{=3\log 2+\frac{15}{2}}b_3 +\underbrace{\(-\frac{3}{2}(2\log 2+1)-4\)}_{=-3\log 2 -\frac{11}{2}}b_5\)\\&
+\sqt\( \underbrace{\(-\frac{3}{2}+\frac{1}{2}\)}_{=-1}c_1 +6 c_3 -6c_5\) \\& +\sqt\( \underbrace{\( \frac{3}{2}-\frac{1}{2}\)}_{ =1} d_1+ \underbrace{\(-\frac{3}{2}-6+\frac{1}{2}-12\)}_{-19}d_3+\underbrace{\(6+6+36\)}_{=48}d_5+\underbrace{\(-6-24\)}_{=-30}d_7\)\\&+\underbrace{\(-\frac{9}{2}-\frac{3}{2}\)}_{=-6}r_3+12r_5-6r_7 +\underbrace{\(\frac{9}{2}+\frac{3}{2}\)}_{=6}s_1-12s_3+6s_5\\& +\underbrace{(3+1)}_{=4}e_3-6e_5 +\underbrace{(-3-1)}_{=-4}f_1  +\underbrace{(9+1)}_{=10}f_3-6f_5
\end{align*}
\normalsize
which concludes the proof of Lemma  \ref{lem:gamma1}.

\qed

\begin{lemma}\label{lem:gamma2} We have %\small
  \begin{align*}
\gamma_2&= \sqt\(-\frac{1}{4}(2\log 2+1)b_1+\(  \log2-\frac{3}{2} \)b_3+4b_5-c_1+2c_3 +d_1-7d_3+6d_5\)\\&\quad+2s_3-2f_3+2f_5 .
\end{align*}
%\normalsize
\end{lemma}
\proof We write  \begin{align*}
\gamma_2&=2\<\Delta_2,h_{3,2}\>=2\<F\xi_{3,1}\xi_{3,2}+\phi_3R_1\xi_{3,1}+\phi_3\xi_{3,1}R_2,h_{3,2}\>\\&
=2\<F\xi_{3,1}\xi_{3,2},h_{3,2}\>+2\<\phi_3 R_1 \xi_{3,2},h_{3,2}\>+2\<\phi_3 \xi_{3,1}R_{2},h_{3,2}\>\\&
=:\gamma_{21}+\gamma_{22}+\gamma_{23}.
\end{align*}

\begin{claim}
   We have
   \begin{align}\label{eq:gamma21}
     \gamma_{21}=-\frac{1}{2\sqt}(2\log 2+1)(b_1-2b_3)-\sqt\(c_1-2c_3\)+\sqt \(d_1-3d_3+2d_5\) .
   \end{align}
\end{claim}
\proof It is a consequence of
\begin{align*}
\gamma_{21}&=2\<F\xi_{3,1}\xi_{3,2},h_{3,2}\>\\&
=2\<\(\frac{1}{4\sqt}(2\log 2+1)\sech +\frac{1}{\sqt}\sech \log \circ \sech -\frac{1}{\sqt}\sech\cdot x\tanh\)(1-2\sechs),-\tanh\sin\>\\&
=-\frac{1}{2\sqt}(2\log 2+1)\int \sech (1-2\sechs)\tanh\sin \\&\quad
-\frac{2}{\sqt}\int \sech(1-2\sechs) \log\circ \sech \tanh\sin\\&\quad
+\sqt\int\sech(1-3\sechs+2\sech^4)  x\sin .
\end{align*}
\qed

\begin{claim}
   We have
   \begin{align}\label{eq:gamma22}
     \gamma_{22}=  -4\sqt (d_3-d_5) +\frac{1}{2}(3s_1-2s_3)-f_1+f_3   .
   \end{align}
\end{claim}
\proof It is  a consequence of
\begin{align*}
\gamma_{22}&= 2\<\phi_3 R_1 \xi_{3,2},h_{3,2}\> \\&
=-2\sqt\<\sech \(2x\tanh \sech^2-\frac{1}{4\sqt}(3-2\sechs)T+\frac{1}{2\sqt}\tanh T'\),\tanh\sin\>\\&
=-4\sqt \int \sech^3(1-\sechs)x\sin +\frac{1}{2}\int (3\sech-2\sech^3)T\tanh\sin -\int \sech(1-\sechs)T'\sin.
\end{align*}

\qed

\begin{claim}
   We have
   \begin{align}\label{eq:gamma23}
     \gamma_{23}=  -2\sqt (b_3-2b_5)-\frac{3}{2}(s_1-2s_3)+f_1-3f_3+2f_5    .
   \end{align}
\end{claim}
\proof It is  a consequence of \small
\begin{align*}&
\gamma_{23} =  2\<\phi_3 \xi_{3,1}R_{2},h_{3,2}\> \\&=-2\sqt\<\sech(1-2\sechs)\(\sechs  +\frac{3}{4\sqt}T-\frac{\tanh}{2\sqt}T'\) ,\tanh\sin\>\\&
=-2\sqt \int (\sech^3-2\sech^5)\tanh\sin
-\frac{3}{2}\int (\sech-2\sech^3)T\tanh\sin
+\int (\sech-3\sech^3+2\sech^5)T'\sin .
\end{align*} \normalsize
\qed

Summing up the quantities computed in the three claims  we get Lemma \ref{lem:gamma2}.
\qed

\begin{lemma}\label{lem:gamma3} We have \small
  \begin{align*}
\gamma_3&= 6\sqt (a_5-2a_7)-6\sqt (d_3-3d_5+2d_7)
 -\frac{7}{\sqt}(p_5-2p_7)+\frac{7}{\sqt}(b_3-2b_5)  .
\end{align*}
\normalsize
\end{lemma}
\proof  We have
\begin{align*}
\gamma_3&=\<(6x\tanh \sechs -\frac{7}{2}\sechs)(\phi_3\xi_{3,1}v_{3,2}),h_{3,1}\>\\&
=6\<x\tanh\sechs \phi_3 \xi_{3,1}\xi_{3,2},h_{3,1}\>-\frac{7}{2}\<\sechs \phi_3 \xi_{3,1 }\xi_{3,2},h_{3,1}\>\\&
=:\gamma_{31}+\gamma_{32}.
\end{align*}

\begin{claim}
   We have
   \begin{align}\label{eq:gamma31}
     \gamma_{31}=  6\sqt (a_5-2a_7)-6\sqt (d_3-3d_5+2d_7)   .
   \end{align}
\end{claim}
\proof It is a consequence of %\small
\begin{align*}
\gamma_{31}&=6\<x\tanh\sechs \phi_3 \xi_{3,1}\xi_{3,2},h_{3,1}\>\\&
=6\sqt\<x\tanh\sech^3 (1-2\sechs),\sechs\cos-\tanh\sin\>\\&
=6\sqt\int \sech^5(1-2\sechs) x\tanh\cos -6\sqt\int \sech^3(1-3\sechs+2\sech^4)x\sin.
\end{align*}% \normalsize
\qed

\begin{claim}
   We have
   \begin{align}\label{eq:gamma31}
     \gamma_{32}= -\frac{7}{\sqt}(p_5-2p_7)+\frac{7}{\sqt}(b_3-2b_5)  .
   \end{align}
\end{claim}
\proof It is a consequence of \small
 \begin{align*}
\gamma_{32}&=-\frac{7}{2}\<\sechs \phi_3 v_{1,3}v_{2,3},h_{3,1}\>\\&
=-\frac{7}{\sqt}\<\sech^3 (1-2\sechs),\sechs\cos-\tanh\sin\>\\&
=-\frac{7}{\sqt}\int \sech^5(1-2\sechs)\cos +\frac{7}{\sqt}\int \sech^3(1-2\sechs)\tanh\sin.
\end{align*} \normalsize
\qed

Summing up the formulas in the Claims we get the proof of Lemma \ref{lem:gamma3}. \qed

\begin{lemma}\label{lem:gamma4} We have
  \begin{align*}
\gamma_4&=   \sqt q_3 -\sqt c_1 + \sqt a_3-\sqt (d_1-d_3)   .
\end{align*}

\end{lemma}
\proof We have
\begin{align*}
\gamma_4&=-2\<E,h_{3,1}\>\\&=-2\<\frac{1}{\sqt} \sech \(\cancel{-\frac{1}{2}\log 2 + \frac{1}{4}}  -\log \sech - x\tanh  \),h_{3,1}\>\\&
=\sqt\<\sech\log\circ \sech,h_{3,1}\>+\sqt\<\sech\cdot x\tanh,h_{3,1}\>\\&
=:\gamma_{42}+\gamma_{43} ,
\end{align*}
 where the canceled term is null by    $\< \phi _3, h_{3,1}\> =0$.  Then the statement follows from
\begin{align*}
\gamma_{42}&=\sqt\<\sech\log\circ\sech,h_{3,1}\>\\&
=\sqt\<\sech\log\circ \sech,\sechs\cos-\tanh\sin\>\\&
=\sqt\int \sech^3 \log\circ\sech  \cos -\sqt\int \sech \log\circ\sech  \tanh\sin\\&
=\sqt q_3 -\sqt c_1
\end{align*}
and
\begin{align*}
\gamma_{43}&=\sqt\<\sech\cdot x\tanh,h_{3,1}\>\\&
=\sqt\<\sech\cdot x\tanh,\sechs\cos-\tanh\sin\>\\&
=\sqt\int \sech^3\cdot x\tanh \cos-\sqt \int \sech (1-\sechs) x\sin\\&
=\sqt a_3-\sqt (d_1-d_3)  ,
\end{align*}
where we remind the reader that our convention is that $\sech ^n \cdot x$ is the product of the function $\sech ^n  (x)$
with the function $x$.

\qed

We now consider the following reduction formulas.
\begin{lemma}\label{lem:red001}
  We have the following relations:
  \begin{align*}
    b_k&
=(k+1)p_{k+2}-kp_k;\\ c_k&
=(k+1)q_{k+2} -kq_k +p_{k+2}-p_k;\\ d_k&
=-ka_k+p_k;\\ e_k&
=s_k+kr_k-(k+1)r_{k+2};\\ f_k&
=-r_k+ks_k.
  \end{align*}
  \end{lemma}
\proof  The formulas follow from the following ones:
\begin{align*}
b_k&=\int \sech^k \tanh \sin \\&
=-k\int\sech^k(1-\sechs)\cos +\int \sech^{k+2}\cos  \\&=
-kp_k+kp_{k+2}+p_{k+2}  ;
\end{align*}
\begin{align*}
c_k&=\int \sech^k \log\circ \sech   \tanh \sin \\&
=-k\int \sech^k (1-\sechs) \log\circ \sech \cos - \int \sech^k (1-\sechs) \cos +\int \sech^{k+2} \log\circ \sech  \cos\\&
=-k(q_k-q_{k+2}) -p_k+p_{k+2} +q_{k+2} ;
\end{align*}
\begin{align*}
d_k&=\int \sech^k \cdot x\sin \\&=
-k\int \sech^k   x \tanh \cos +\int \sech^k\cos\\&
=-ka_k+p_k;
\end{align*}\begin{align*}
e_k&=\int \sech^k T' \tanh \cos\\&=\int \sech^k T\tanh \sin  -\int \sech^{k+2}T\cos + k\int \sech^k(1-\sechs)T\cos\\&
=s_k-r_{k+2} +k (r_k-r_{k+2}) ;
\end{align*}
\begin{align*}
f_k&=\int \sech^k T'\sin\\&
=-\int \sech^k T\cos + k\int \sech^k T\tanh \sin .
\end{align*}
\qed

Thanks to Lemma \ref{lem:red001} we can eliminate some of the variables.
\begin{lemma}\label{lem:redform} We have the following formulas:
   \begin{align}\label{eq:gam1} \gamma _1& = \sqt\( 2 p_1+\(-9\log 2-\frac{97}{2}\)p_3+ \(24\log 2+110\)p_5+\(-15\log 2 -\frac{127}{2}\)p_7\)\\&\quad
+\sqt \(  q_1 -19 q_3+48 q_5 -30 q_7-a_1+ 56 a_3 -222 a_5+180 a_7\) \nonumber \\&\quad
+4r_1-4r_3 -28 r_5 +30r_7 +2 s_1 +22 s_3 -30 s_5 \nonumber
    ;\\
%   \end{align}
%    \begin{align}
    \label{eq:gam2} \gamma _2& = \sqt\(\(\frac{1}{2}\log2 +\frac{9}{4}\)p_1+\(-4\log2-6\)p_3+ \(4\log2 -18\) p_5+24p_7 \) \\& \quad+ \sqt \( q_1-8q_3+8q_5-a_1+21a_3-30a_5\)   \nonumber \\&\quad+
2r_3-2r_5-4s_3+10s_5 \nonumber
    ;\\
%   \end{align}
%   \begin{align}
\label{eq:gam3} \gamma _3& =  \sqt\(-\frac{33}{2}p_3+\frac{127}{2}p_5-47p_7+18a_3-84a_5+72a_7\)
    ;\\
%   \end{align}
%   \begin{align}
   \label{eq:gam4} \gamma _4& =  \sqt\(q_1-q_3+a_1-2a_3\) .
   \end{align}

\end{lemma}

\proof Starting from the simplest formulas we have
\begin{align*}
\gamma_4&
=\sqt q_3 -\sqt c_1+\sqt a_3-\sqt d_1+\sqt d_3\\&
=\sqt q_3 -\sqt (2q_3-q_1+\cancel{p_3}\cancel{-p_1})+\sqt a_3-\sqt(-a_1+\cancel{p_1})+\sqt (-3a_3+\cancel{p_3})\\&
=\sqt q_1-\sqt q_3+\sqt a_1-2\sqt a_3\\&
=\sqt\(q_1-q_3+a_1-2a_3\),
\end{align*}
which yields \eqref{eq:gam4}. Then we consider
\small \begin{align*}
&\gamma_3 =
6\sqt (a_5-2a_7)-6\sqt (d_3-3d_5+2d_7)
-\frac{7}{\sqt}(p_5-2p_7)+\frac{7}{\sqt}(b_3-2b_5)\\&
=\sqt \( 6 a_5-12a_7- 6d_3+18d_5-12d_7
-\frac{7}{2}p_5+7p_7+\frac{7}{2}b_3-7b_5\)\\&
=\sqt\(6a_5-12a_7-6 (-3a_3+p_3)+18(-5a_5+p_5)-12(-7a_7+p_7)
-\frac{7}{2}p_5\right . \\&\quad\quad \left . +7p_7+\frac{7}{2}(4p_5-3p_3)-7(6p_7-5p_5)\)\\&
=\sqt\(\underbrace{(-6-\frac{21}{2})}_{=-\frac{33}{2}}p_3+(\underbrace{18-\frac{7}{2}+14+35}_{=\frac{127}{2}})p_5+\underbrace{(-12+7-42)}_{=-47}p_7+18a_3+\underbrace{(6-90)}_{=-84}a_5+\underbrace{(-12+84)}_{=72}a_7\)\\&
=\sqt\(-\frac{33}{2}p_3+\frac{127}{2}p_5-47p_7+18a_3-84a_5+72a_7\) .
\end{align*}
\normalsize
Using Lemma \ref{lem:gamma2} we have   \small
\begin{align*}
&\gamma_2 = \sqt\(-\frac{1}{4}(2\log 2+1)b_1+\(  \log2-\frac{3}{2} \)b_3+4b_5-c_1+2c_3 +d_1+(-3-4)d_3+(2+4)d_5\)\\& \quad+2s_3-2f_3+2f_5
\\&  = \sqt\(-\frac{1}{4}(2\log 2+1)(2p_3-p_1)+\(\log2-\frac{3}{2}\)(4p_5-3p_3)+4(6p_7-5p_5)-(2q_3-q_1+p_3-p_1)\)\\&\quad+\sqt\(2(4q_5-3q_3+p_5-p_3) +(-a_1+p_1)-7(-3a_3+p_3)+6(-5a_5+p_5)\)\\&\quad+2s_3-2(-r_3+3s_3)+2(-r_5+5s_5).
\end{align*}
\normalsize
Collecting together similar terms we have the following, which yields \eqref{eq:gam3},  \small
\begin{align*}
&\gamma_2 = \sqt\(\underbrace{\(\frac{1}{4}(2\log 2+1)+1+1\)}_{=\frac{1}{2}\log2 +\frac{9}{4}}p_1+\underbrace{\(-\frac{1}{2}(2\log 2+1)-3\(\log2-\frac{3}{2}\)-1-2-7\)}_{=-4\log2-6}p_3\)\\&\quad
+\sqt\( \underbrace{\(4\(\log2-\frac{3}{2}\)-20+2+6\)}_{=4\log2 -18} p_5+24p_7+q_1+(-2-6)q_3+8q_5-a_1+21a_3-30a_5\)\\&\quad+
2r_3-2r_5+(2-6)s_3+10s_5 .
\end{align*}
\normalsize
Using Lemma \ref{lem:gamma1} and Lemma \ref{lem:red001} we have   \small
\begin{align*}
\gamma_1 = & \sqt\(  \(-3\log2-\frac{15}{2}\)p_5+\(3\log2+\frac{11}{2}\)p_7+ q_3 -6q_5+6q_7 -a_3 + 18 a_5-30 a_7\)\\&\quad
+\sqt\( \(3\log 2+\frac{15}{2}\)(-3p_3+4p_5) +\(-3\log 2 -\frac{11}{2}\)(-5p_5+6p_7)\)\\&\quad
+\sqt\( -(-p_1+p_3-q_1+2q_3) +6 (-p_3+p_5-3q_3+4q_5) -6(-p_5+p_7-5q_5+6q_7)\) \\&\quad
+\sqt\((p_1-a_1) -19 (p_3-3a_3)+ 48 (p_5-5a_5)-30 (p_7-7 a_7)\)\\&\quad
-6r_3+12 r_5-6r_7 +6 s_1-12s_3+6s_5+ 4 (3r_3-4r_5+s_3)-6(5r_5-6r_7+s_5) \\&\quad-4 (-r_1+s_1)  +10(-r_3+3s_3)-6(-r_5+5s_5).
\end{align*}
\normalsize
Collecting similar terms we have the following, which yields \eqref{eq:gam1},
\small
\begin{align*}
\gamma_1 =&  \sqt\( \underbrace{\(1+1\)}_{=2}p_1+\underbrace{\(-3\(3\log 2+\frac{15}{2}\)-1-6-19\)}_{=-9\log 2-\frac{97}{2}}p_3\) \\&
+\sqt\( \underbrace{\( \(-3\log2-\frac{15}{2}\)+4\(3\log 2 + \frac{15}{2}\)-5\(-3\log 2 -\frac{11}{2}\)+6+6+48\)}_{=24\log 2+110}p_5\)\\&  +\sqt\(\underbrace{\(\(3\log2+\frac{11}{2}\) +6\(-3\log 2 -\frac{11}{2}\)-6-30\)}_{=-15\log 2 -\frac{127}{2}}p_7\)\\&
+\sqt \(  q_1 + \(1-2-18\) q_3+ \(-6+24+30\) q_5+ \(6-36\) q_7 \)\\& +\sqt \(  -a_1+ (-1+57) a_3+ (18-240) a_5+  (-30+210) a_7\)\\&
+4r_1+ (-6+12-10) r_3+ (12-16-30+6) r_5+ (36-6) r_7 \\& + (6-4) s_1+ (-12+4+30) s_3+(\cancel{6-6}-30)s_5.
\end{align*}
\normalsize
\qed

\textit{Proof of Proposition \ref{lem:1st}.} Summing up the formulas in Lemma \ref{lem:redform}, we obtain  \small
\begin{align*}
\Gamma =&\gamma_1+\gamma_2+\gamma_3+\gamma_4\\
=&\sqt\( 2 p_1+\(-9 \log 2-\frac{97}{2}\)p_3+ \(24\log 2+110\)p_5+\(-15\log 2 -\frac{127}{2}\)p_7\)\\&
+\sqt \(  q_1 -19 q_3+48 q_5 -30 q_7-a_1+ 56 a_3 -222 a_5+180 a_7\)\\&
+4r_1-4r_3 -28r_5 +30r_7 +2 s_1 +22 s_3 -30 s_5\\&
+\sqt\(\(\frac{1}{2}\log2 +\frac{9}{4}\)p_1+\(-4\log2-6\)p_3+ \(4\log2 -18\) p_5+24p_7+q_1-8q_3+8q_5 \)\\&
+\sqt\( -a_1+21a_3-30a_5\) +
2r_3-2r_5-4s_3+10s_5\\&\quad
+\sqt\(-\frac{33}{2}p_3+\frac{127}{2}p_5-47p_7+18a_3-84a_5+72a_7\)\\&
+\sqt\(q_1-q_3+a_1-2a_3\)
\end{align*}
\normalsize
where collecting together similar terms we have the following which yields
Proposition \ref{lem:1st}
 \small
\begin{align*}
\Gamma
=&\sqt\(\underbrace{\(\frac{1}{2}\log 2 +2+\frac{9}{4}\)}_{=\frac{1}{2}\log 2 + \frac{17}{4}}p_1 +  \underbrace{\(\(-9-4\)\log 2 -\frac{97}{2}-6-\frac{33}{2}\)}_{=-13\log 2 -71}p_3\)\\&\quad
+\sqt\( \underbrace{\(  \(24+4\)\log 2 + 110-18+\frac{127}{2}\)}_{=28\log 2 +\frac{311}{2}}p_5 + \underbrace{\( -15\log 2 -\frac{127}{2} +24-47\)}_{=-15\log 2 -\frac{173}{2}}p_7\)\\&\quad
+\sqt\(\underbrace{(1+1+1)}_{=3}q_1 + \underbrace{(-19-8-1)}_{=-28}q_3 + \underbrace{(48+8)}_{=56}q_5 -30q_7\)\\&\quad
+\sqt\(\underbrace{(-1-1+1)}_{=-1}a_1 + \underbrace{(56+21+18-2)}_{=93}a_3+\underbrace{(-222-30-84)}_{=-336}a_5+\underbrace{(180+72)}_{252}a_7\)\\&\quad
+4r_1+\underbrace{(-4+2)}_{=-2}r_3+\underbrace{(-28-2)}_{=-30}r_5 +30r_7 +2s_1 + \underbrace{(22-4)}_{=18}s_3 +\underbrace{(-30+10)}_{=-20}s_5.
\end{align*}
\normalsize
\qed

\section{Reductions and cancellations} \label{sec:cancel}

Now we want to perform further reductions and express $\Gamma$ in terms of $p_1$, $q_1$, $a_1$, $r_1$ and $s_1$.

\begin{lemma} \label{lem:redp}
  We have the formulas \begin{align*}
p_3&=  p_1,\\
p_5&= \frac{5}{6}p_1,\\
p_7&= \frac{13}{18}p_1.
\end{align*}
\end{lemma}
\proof We have \begin{align*}
p_k&=\int \sech^k (\sin)'
=k\int \sech^k\tan  (-\cos)'
=k\int (\sech^k \tan)' \cos \\&
=-k^2\int \sech^k(1-\sechs)\cos + k\int \sech^{k+2}\cos\\&
=-k^2 p_k + (k^2+k)p_{k+2} .
\end{align*}
Thus,
\begin{align*}
p_{k+2}=\frac{1+k^2}{k(k+1)}p_k
\end{align*}
and
\begin{align*}
p_3&=\frac{1+1}{1\cdot(1+1)}p_1=p_1,\\
p_5&=\frac{1+3^2}{3\cdot 4}p_3=\frac{10}{12}p_1=\frac{5}{6}p_1,\\
p_7&=\frac{1+5^2}{5\cdot 6}p_5=\frac{26}{30}\frac{5}{6}p_1=\frac{13}{18}p_1.
\end{align*}
\qed

\begin{lemma} \label{lem:redq}
  We have the formulas \begin{align*}
q_3&= q_1-\frac{1}{2}p_1,\\
q_5&=  \frac{5}{6}q_1-\frac{29}{72}p_1,\\
q_7&
=\frac{13}{18}q_1-\frac{121}{360}p_1 .
\end{align*}
\end{lemma}
\proof  Using $(\log \sech)'=\frac{1}{\sech} (-\sech\tanh)=-\tanh$, we have
\begin{align*}
q_k&=\int \sech^k \log\circ\sech (\sin)'
=k\int \sech^k \log\circ\sech \tanh \sin +\int \sech^k \tanh \sin\\&
=kc_k+b_k\\&
=k\((k+1)q_{k+2}-kq_k+p_{k+2}-p_k\) +(k+1)p_{k+2}-kp_k\\&
=k(k+1)q_{k+2}-k^2 q_k + (2k+1)p_{k+2} -2k p_k.
\end{align*}
Thus,
\begin{align*}
q_{k+2}=\frac{1}{k(k+1)}\((1+k^2)q_k -(2k+1)p_{k+2} + 2k p_k\).
\end{align*}
In particular,
\begin{align*}
q_3&=\frac{1}{2}\((1+1)q_1-(2+1)p_3 + 2p_1\)=q_1-\frac{1}{2}p_1,\\
q_5&=\frac{1}{3\cdot 4}\((1+9)q_3 -(2\cdot 3+1)p_5 +2\cdot 3 p_3\)=\frac{5}{6}q_3-\frac{7}{12}p_5 +\frac{1}{2}p_3\\&
=\frac{5}{6}\(q_1-\frac{1}{2}p_1\)-\frac{7}{12}\frac{5}{6}p_1+\frac{1}{2}p_1 =\frac{5}{6}q_1-\frac{29}{72}p_1,\\
q_7&=\frac{1}{5\cdot 6}\((1+5^2)q_5 -(10+1)p_7 +2\cdot 5 p_5\)\\&
=\frac{13}{15}q_5 -\frac{11}{30}p_7 +\frac{1}{3}p_5=\frac{13}{15}\(\frac{5}{6}q_1-\frac{29}{72}p_1\)-\frac{11}{30}\frac{13}{18}p_1+\frac{1}{3}\frac{5}{6}p_1\\&
=\frac{13}{18}q_1-\frac{121}{360}p_1 .
\end{align*}
\qed

\begin{lemma} \label{lem:redsr}
  We have the formulas \begin{align*}
r_3&= -r_1+s_1+\sqt p_1,\\
s_3&
=\frac{1}{3}s_1-r_1+\frac{7}{9}\sqt p_1, \\
r_5&
=-r_1 + \frac{2}{3}s_1 +\frac{37}{36}\sqt p_1,\\
s_5&
=-\frac{2}{5}r_1 +\frac{1}{15}s_1 +\frac{13}{36}\sqt p_1, \\
r_7&= -\frac{13}{15}r_1 + \frac{23}{45}s_1+\frac{83\sqt}{90} p_1.
\end{align*}\end{lemma}

\proof By $T''=2T-2\sqt \sechs$ and $e_k=kr_k-(k+1)r_{k+2}+s_k$.
\begin{align*}
r_k&=\int \sech^k T (\sin)'\\&
=-\int (\sech^k T)' \sin\\&
=k\int \sech^k \tanh T \sin -\int \sech^k T' \sin\\&
=ks_k - \int \sech^k T' (-\cos)'\\&
=ks_k -\int (\sech^k T')'\cos\\&
=ks_k +k\int \sech^k\tanh T'\cos -\int \sech^k(2T-2\sqt\sechs)\cos\\&
=ks_k + ke_k -2r_k+2\sqt p_{k+2}\\&
=ks_k + k(kr_k-(k+1)r_{k+2}+s_k) -2r_k+2\sqt p_{k+2}\\&
=-k(k+1)r_{k+2}+(k^2-2) r_k+2ks_k+2\sqt p_{k+2}.
\end{align*}
Thus,
\begin{align*}
r_{k+2}=\frac{1}{k(k+1)}\((k^2-3) r_k+2ks_k+2\sqt p_{k+2}\) .
\end{align*}
 We have
\begin{align*}
s_k&=\int \sech^k T \tanh (-\cos)'\\&
=\int (\sech^k T \tanh)' \cos\\&
=-k\int \sech^k(1-\sechs)T\cos + \int \sech^{k+2}T \cos +\int \sech^k T'\tanh {\cos} \\&
=-kr_k + (k+1)r_{k+2} -\int (\sech^k T' \tanh)' \sin\\&
=-kr_k + (k+1)r_{k+2} +k\int \sech^k (1-\sechs)T'  \sin \\& -\int \sech^{k+2} T'\sin -\int \sech^k(2T-2\sqt \sechs)\tanh \sin \\&
=-kr_k + (k+1)r_{k+2} +k(f_k-f_{k+2}) -f_{k+2} -2s_k +2\sqt b_{k+2} \\&
=-kr_k + (k+1)r_{k+2} +k(-r_k+ks_k)-(k+1)(-r_{k+2}+(k+2)s_{k+2})\\& \quad  -2s_k +2\sqt ((k+3)p_{k+4}-(k+2)p_{k+2}) .
\end{align*}
Thus,
\begin{align*}
s_{k+2}=\frac{1}{(k+1)(k+2)}\( (k^2-3) s_k +2(k+1)r_{k+2} -2k r_k +2\sqt (k+3) p_{k+4}-2\sqt (k+2) p_{k+2}\) .
\end{align*}
Hence
\small
\begin{align*}
r_3&=\frac{1}{2}\((1-3)r_1+2s_1+2\sqt p_3\) =-r_1+s_1+\sqt p_1 ,\\
s_3&=\frac{1}{2\cdot 3} \((1-3)s_1 +2\cdot 2 r_3-2r_1 +2\sqt(1+3)p_5-2\sqt (1+2)p_3 \)\\&
=-\frac{1}{3}s_1 +\frac{2}{3}(-r_1+s_1+\sqt p_1) -\frac{1}{3}r_1 +\frac{4\sqt}{3}\cdot \frac{5}{6}p_1 -\sqt p_1\\&
=\frac{1}{3}s_1-r_1+\frac{7}{9}\sqt p_1 ,\end{align*}
\begin{align*}r_5&=\frac{1}{12}\((9-3)r_3 +6s_3 \)
=\frac{1}{3\cdot 4}\((9-3)r_3+2\cdot 3 s_3 +2\sqt p_5\)
=\frac{1}{12}\(6r_3+6s_3+2\sqt p_5\)\\&
=\frac{1}{2}\(-r_1+s_1+\sqt p_1\)+\frac{1}{2}\(\frac{1}{3}s_1-r_1+\frac{7}{9}\sqt p_1\)+\frac{\sqt}{6}\cdot \frac{5}{6}p_1\\&
=-r_1 + \frac{2}{3}s_1 +\frac{37}{36}\sqt p_1,\\
s_5&=\frac{1}{4\cdot 5}\(6s_3 +8 r_5 -6r_3 +12\sqt p_7-10\sqt p_5\)\\&=
\frac{3}{10}\(\frac{1}{3}s_1-r_1+\frac{7}{9}\sqt p_1\)+\frac{2}{5}\(-r_1 + \frac{2}{3}s_1 +\frac{37}{36}\sqt p_1\)-\frac{3}{10}\(-r_1+s_1+\sqt p_1\) \\& \quad+\frac{3}{5}\sqt \frac{13}{18}p_1 -\frac{1}{2}\sqt \frac{5}{6}p_1  \\&
=-\frac{2}{5}r_1 +\frac{1}{15}s_1 +\frac{13}{36}\sqt p_1 \end{align*} and
\begin{align*}r_7&=\frac{1}{5\cdot 6}\((25-3)r_5+10s_5+2\sqt p_7\)\\&
=\frac{11}{15}\(-r_1 + \frac{2}{3}s_1 +\frac{37}{36}\sqt p_1\) +\frac{1}{3}\(-\frac{2}{5}r_1 +\frac{1}{15}s_1 +\frac{13}{36}\sqt p_1\) +\frac{\sqt}{15}\frac{13}{18}p_1\\&
=-\frac{13}{15}r_1 + \frac{23}{45}s_1+\frac{83\sqt}{90} p_1 .
\end{align*}
\normalsize

\qed

\textit{Proof of Theorem \ref{thm:fgr}.} We substitute in the formula of $\Gamma$ in Proposition \ref{lem:1st} the formulas in
Lemmas \ref{lem:redp}, \ref{lem:redq} and \ref{lem:redsr}. We obtain \small
\begin{align*}
\Gamma
=&\sqt\(\(\frac{1}{2}\log 2 + \frac{17}{4}\)p_1 -  \( 13\log 2 +71\)p_1+ \(28\log 2 +\frac{311}{2}\)\frac{5}{6}p_1 - \( 15\log 2 +\frac{173}{2}\)\frac{13}{18}p_1\)\\&\quad
+\sqt\(3q_1 -28\(q_1-\frac{1}{2}p_1\) +56 \(\frac{5}{6}q_1-\frac{29}{72}p_1\) -30\(\frac{13}{18}q_1-\frac{121}{360}p_1\)\)\\&\quad+\sqt\(- a_1 + 93 \(\frac{1}{3}a_1+\frac{1}{3}p_1\)-336 \(\frac{1}{6}a_1 +\frac{1}{5}p_1\) +252 \(\frac{13}{126} a_1 + \frac{83}{630}p_1\)\)\\&\quad
+4r_1 - 2 \(-r_1+s_1+\sqt p_1\) - 30 \(-r_1 + \frac{2}{3}s_1 +\frac{37}{36}\sqt p_1\) +30\(-\frac{13}{15}r_1 + \frac{23}{45}s_1+\frac{83\sqt}{90} p_1\) \\&\quad+2s_1 +  18 \(\frac{1}{3}s_1-r_1+\frac{7\sqt}{9}p_1\) -20 \(-\frac{2}{5}r_1 +\frac{1}{15}s_1 +\frac{13\sqt}{36}p_1\) .
\end{align*}\normalsize
Collecting together similar terms and cancelling the null ones, we obtain
\begin{align*}
&\Gamma
=\sqt  \alpha p_1    + \sqt\cancel{ \(3-28+56\frac{5}{6}-30 \frac{13}{18}\)} q_1 \\&\quad+ \cancel{\(4+2+\frac{63}{2}-\frac{30 \cdot 13}{15}-18+8\)}r_1 + \cancel{\(-2-\frac{63\cdot 2}{2\cdot 3} +30\cdot \frac{23}{45} +2+\frac{18}{3}-20\cdot \frac{1}{15}\)}s_1\\&\quad +\sqt \cancel{\(-1+\frac{93}{3} -\frac{336}{6} +\frac{252\cdot 13}{126}\)}a_1
\end{align*}
where \small
\begin{align*}
 & \alpha =  \log 2\cancel{\(\frac{1}{2}-13 + 28\frac{5}{6}-15 \frac{13}{18}\)} + \\&   \frac{17}{4}-71 +\frac{311}{2}\frac{5}{6} -\frac{173\cdot 13}{2\cdot 18} +14-56\cdot \frac{29}{72} +30 \frac{121}{360} \\& +\frac{93}{3}-\frac{336}{5} +252\frac{83}{630} -2-30\frac{37}{36}+30\frac{83}{90} +18\frac{7}{9} -20\frac{13}{36}  =\frac{1}{2}
\end{align*}
\normalsize
Hence $\Gamma = \frac{1}{\sqrt{2}}  p_1 $. Finally, $p_1=\pi \sech (\pi/2)$ by an application of the Residue Theorem.

\qed

We remark that we do not have a conceptual justification of why the computations  give such a simple formula for $\Gamma$ which we would not expect from the outset.

\section*{Acknowledgments}
C. was supported   by the Prin 2020 project \textit{Hamiltonian and Dispersive PDEs} N. 2020XB3EFL.
M.  was supported by the JSPS KAKENHI Grant Number 19K03579, G19KK0066A, 23H01079 and 24K06792.

Department of Mathematics and Geosciences,  University
of Trieste, via Valerio  12/1  Trieste, 34127  Italy.
{\it E-mail Address}: {\tt scuccagna@units.it}

Department of Mathematics and Informatics,
Graduate School of Science,
Chiba University,
Chiba 263-8522, Japan.
{\it E-mail Address}: {\tt maeda@math.s.chiba-u.ac.jp}
\end{document}